\documentclass[12pt]{article}
\usepackage{amssymb,amsmath}
\usepackage{graphicx}
\usepackage{amsfonts}
\usepackage{float}
\usepackage{flafter}
\usepackage{color}

\labelwidth\leftmargini\advance\labelwidth-\labelsep

\def\R{\relax\ifmmode I\!\!R\else$I\!\!R$\fi}

\def\Z{\relax\ifmmode Z\!\!\!Z\else$Z\!\!\!Z$\fi}

\def\C{\relax\ifmmode C\!\!\!\!I\else$C\!\!\!\!I$\fi}

\def\K{\relax\ifmmode I\!\!K\else$I\!\!K$\fi}

\def\N{\relax\ifmmode I\!\!N\else$I\!\!N$\fi}

\newcounter{defcounter}[section]
\newenvironment{definition}%
{\vspace{0.1cm}\begin{sloppypar}\noindent\stepcounter{defcounter}{\bfseries
Definition
      \thesection.\thedefcounter}}%
{\end{sloppypar}\vspace{0.1cm}}

\newtheorem{lemma}{Lemma}[section]
\newtheorem{theorem}{Theorem}[section]

\newtheorem{proposition}{Proposition}[section]

\newcommand{\proof}{{\bf Proof.} }

\newcommand{\qed}{\hfill $\square$}

\begin{document}
\thispagestyle{empty}
\begin{center}
{\Large {\bf A new family of expansions of real numbers}}
\end{center}
\begin{center}J\"org Neunh\"auserer\\
Technical University of Braunschweig \\
joerg.neunhaeuserer@web.de
\end{center}
\begin{center}
\begin{abstract}
For $\alpha>1$ we represent a real number in $(0,1]$ in the form
\[ \sum_{i=1}^{\infty}(\alpha-1)^{i-1}\alpha^{-(d_{1}+\dots+d_{i})}\]
with $d_{i}\in\mathbb{N}$. We discuss ergodic theoretical and dimension theoretical aspects of this expansion. Furthermore we study their base-change-transformation.  
\end{abstract}
\end{center}
\section{Introduction}
Expansion of real numbers, like classical $b$-adic expansions and continued fraction expansions, are a central topic in the metric theory of numbers. We like to remind the reader of some classical results. Borel \cite{[BO]} proved that almost all real numbers are normal to all bases $b$, in the $b$-adic expansions all digits appear with same frequency. It was Khinchin \cite{[KI]} who proved, that for almost all real numbers the asymptotic geometric mean of the digits of the continued fraction expansion is a universal constant, now called Khinchin's constant. Such result are today a part of ergodic theory. Hausdorff \cite{[HA]} introduced his dimension to distinguish the size of real numbers with digits in the $b$-adic expansion from a proper subset of all digits. Furthermore Besicovitch \cite{[BE]} and Eggleston \cite{[EG]} found the Hausdorff dimension of the set of real numbers with given frequencies of digits in the $b$-adic expansion. The dimension theory of continued fraction expansions goes back to Jarnick \cite{[JA]}, who shows that the set of continued fractions with bounded digits has Hausdorff dimension one, although this set has Lebesgue measure zero.\\
Beside  $b$-adic expansions and continued fraction expansions many other expansions of real number are studied. We have ergodic and dimension theoretical results on $\beta$-expansions \cite{[RE],[BA]}, Cantor series expansions  \cite{[ER],[KIF],[DY]}, Engel-expansions  \cite{[GA],[WU],[FANG]}, Lüroth-expansions  \cite{[JAG],[BAR]}, multiple-base expansions \cite{[NE3],[NE3],[LI],[KO]} and others.   \\
In this paper we introduce a new family of expansions of real numbers which is given essentially by the expression in the abstract. We call this expansion $\alpha$-expansion and describe this expansion in detail in the next section. In section 4 we consider ergodic aspects of $\alpha$-expansions and proof an analogon of Borels and of Khinchins theorem. In the next section we study dimension theoretical aspects of $\alpha$-expansions. We proof results similar to the theorems of Hausdorff and Jarnik and similar to the theorems of Besicovitch and Eggleston. In the last section we discuss the base-change-transformation for $\alpha$-expansions. It turns out that the transformation is strictly monoton and continuous with an infinit self-similar graph of Hausdorff dimension one.  Furthermore generically the upper derivative of this transformation is infinite and the lower derivative is zero. 
\section{The expansions} 
Let $\alpha>1$ be a real number. For an infinite sequence $(d_{i})\in\mathbb{N}^{\mathbb{N}}$ we set
\[ \langle d_{i}\rangle_{\alpha}=\langle d_{1},d_{2},\dots\rangle_{\alpha}=\sum_{i=1}^{\infty}(\alpha-1)^{i-1}\alpha^{-(d_{1}+\dots+d_{i})}.\]
Let $\pi:\mathbb{N}^{\mathbb{N}}\to (0,1]$ by given by $\pi((d_{i}))= \langle d_{i}\rangle_{\alpha}$. Since 
\[ \sum_{i=1}^{\infty}(\alpha-1)^{i-1}\alpha^{-(d_{1}+\dots+d_{i})}\le \sum_{i=1}^{\infty}(\alpha-1)^{i-1}\alpha^{-i}=1\]
the map is well defined. We will use an other discription of the map $\pi$. For $i\in\mathbb{N}$ let $T_{i}:[0,1]\to[0,1]$ be given by
\[ T_{i}(x)=\frac{\alpha-1}{\alpha^{i}}x+\frac{1}{\alpha^{i}}.\]
We have:
\begin{lemma}
 \[\pi((d_{i}))=\lim_{n\to\infty}T_{d_{n}}\circ \dots\circ T_{d_{1}}(1).\]
\end{lemma}
\proof By induction we get
\[  T_{d_{n}}\circ \dots\circ T_{d_{1}}(x)=\frac{(\alpha-1)^{n}}{ \alpha^{d_{1}+\dots+d_{n}}}x+ \sum_{i=1}^{n}(\alpha-1)^{i-1}\alpha^{-(d_{1}+\dots+d_{i})}.\]                  The result follows since the maps $T_{i}$ are contractions for all $i\in\mathbb{N}$.
\qed \\\\
Using this lemma we easily obtain:
\begin{proposition}
$\pi$ is a bijection.
\end{proposition}  
\proof
We have
\[\bigcup_{i=1}^{\infty}T_{i}((0,1])=\bigcup_{i=1}^{\infty}\left(\frac{1}{\alpha^{i}},\frac{1}{\alpha^{i-1}}\right]=(0,1]\]
hence $\pi$ is surjective. Moreover
\[ T^{i}((0,1])\cap T^{j}((0,1])=\left(\frac{1}{\alpha^{i}},\frac{1}{\alpha^{i-1}}\right]\cap\left(\frac{1}{\alpha^{j}},\frac{1}{\alpha^{j-1}}\right]=\emptyset \] 
it $i\not=j$, hence $\pi$ is injective. \qed~\\~\\
By this proposition we have an unique expansion of the form $x= \langle d_{i}\rangle_{\alpha}$ with $(d_{i})\in\mathbb{N}^{\mathbb{N}}$ for every $x\in(0,1]$. Thus we may define: 
\begin{definition}
For $\alpha>1$ the sequence $(d_{i})\in\mathbb{N}^{\mathbb{N}}$ with $x= \langle d_{i}\rangle_{\alpha}$ is the $\alpha$-expansion of $x$. 
\end{definition}
\section{Ergodic theoretical aspects}
We use here a dynamical  approach to proof results on the distribution of digits of $\alpha$-expansions. The reader who is not familiar with basic ergodic theory should consider our essential \cite{[NE1]} or \cite{[EIN]}.\\
Let $\sigma:\mathbb{N}^{\mathbb{N}}\to\mathbb{N}^{\mathbb{N}}$ given by $\sigma(d_{k})=d_{k+1}$ be the shift map. We define the corresponding map on $T:(0,1]\to(0,1]$, with respect to the coding map $\pi$, by 
\[ T(x)=T^{-1}_{i}(x) \mbox{ for }x\in T_{i}((0,1])=\left(\frac{1}{\alpha^{i}},\frac{1}{\alpha^{i-1}}\right]\]
for $i\in\mathbb{N}$. By lemma 2.1 we have the conjugation
\[ T\circ\pi =\pi\circ\sigma \]
and obtain
\begin{lemma}
For $x= \langle d_{i}\rangle_{\alpha}\in(0,1]$ we obviously have $d_{j}=d\in\mathbb{N}$, if and only if 
\[ T^{j-1}(x)\in  T_{d}((0,1])=\left(\frac{1}{\alpha^{d}},\frac{1}{\alpha^{d-1}}\right].\] 
\end{lemma} 
\proof Obviously $x= \langle d_{i}\rangle_{\alpha}\in  T_{d}((0,1])$, if and only if $d_{1}=d$. Hence
\[ T^{j-1}(x)= T^{j-1}(\langle d_{i}\rangle_{\alpha})=\langle \sigma^{j-1}(d_{i})\rangle_{\alpha}\in  T_{d}((0,1]),\]
if and only if $d_{j}=d$. \qed~\\~\\ 
To order to apply Birkhoff's ergodic theorem we state:
\begin{proposition}
The Lebesgue measure $\mathfrak{L}$ is ergodic with respect to $T$.
\end{proposition}
\proof
For an open interval $(a,b)\subseteq [0,1]$ we have
\[ \mathfrak{L}(T^{-1}((a,b)))=\mathfrak{L}\left(\bigcup_{i=1}^{\infty}(\frac{\alpha-1}{\alpha^{i}}a+\frac{1}{\alpha^{i}}, \frac{\alpha-1}{\alpha^{i}}b+\frac{1}{\alpha^{i}} \right)\]
\[ = \sum_{i=1}^{\infty}\frac{\alpha-1}{\alpha^{i}}\mathfrak{L}\left((a,b))\right)=b-a=\mathfrak{L}((a,b)).\]
Hence $\mathfrak{L}(T^{-1}(B))=\mathfrak{L}(B)$ for all Borel sets $B\subseteq (0,1]$, which means that $\mathfrak{L}$ is invariant under $T$. We could now directly proof that $\mathfrak{L}$ is ergodic with respect to $T$ using Lebesgue density theorem, but this is not necessary. It is well known that for piecewise smooth expanding interval maps the absolutely continuous invariant measure is in fact ergodic, see \cite{[LY]}.  
\qed~\\~\\
Now we consider the frequency of a digit $d\in\mathbb{N}$ in the $\alpha$-expansion of $x=\langle d_{i}\rangle_{\alpha}\in(0,1]$,
\[ \mathfrak{f}_{d,\alpha}(x):=\lim_{n\to\infty}\frac{\sharp\{j\in\{1,\dots,n\}|d_{j}=d\}}{n}.\]  
Generically this frequencies are geometrically distributed:
\begin{theorem}
Let $\alpha>1$. For almost all $x\in(0,1]$ we have $\mathfrak{f}_{d,\alpha}(x)=(\alpha-1)/\alpha^{d}$ for all $d\in\mathbb{N}$. 
\end{theorem} 
\proof
Let $\chi_{A}$ be the characteristic function of $A\subseteq \mathbb{R}$, that is $\chi_{A}(x)=1$ for $x\in A$ and $\chi_{A}(x)=0$ else. For $x=\langle d_{i}\rangle_{\alpha}\in(0,1]$ we have by lemma 3.1.
\[ \sharp\{j\in\{1,\dots,n\}|d_{j}=d\}=\sum_{i=0}^{n-1}\chi_{\left(\frac{1}{\alpha^{d}},\frac{1}{\alpha^{d-1}}\right]}(T^{i}(x)),\]
for all $d\in\mathbb{N}$. Fix $d\in\mathbb{N}$. By proposition 3.1 and Birkhoffs ergodic theorem we have for almost all $x=\langle d_{i}\rangle_{\alpha}\in(0,1]$
\[ \mathfrak{f}_{d,\alpha}(x)=\lim_{n\to\infty}\frac{\sharp\{j\in\{1,\dots,n\}|d_{j}=d\}}{n}=\lim_{n\to\infty}\frac{1}{n}\sum_{i=0}^{n-1}\chi_{\left(\frac{1}{\alpha^{d}},\frac{1}{\alpha^{d-1}}\right]}(T^{i}(x))\]
\[ =\int_{0}^{1}\chi_{\left(\frac{1}{\alpha^{d}},\frac{1}{\alpha^{d-1}}\right]}(x)dx=(\alpha-1)/\alpha^{d}.\]
The result follows since the countable intersection of sets with full measure has full measure.
\qed\\\\
From this theorem it follows that the set of numbers with a bounded sequence of digits in their $\alpha$-expansion has Lebesgue measure zero. We will see in the next section that the Hausdorff dimension of this set is one.~\\
In our second theorem on the distribution of digits in $\alpha$-expansion we consider the asymptotic arithmetic and geometric mean of the digits.
\begin{theorem} 
Let $\alpha>1$. For almost all $x=\langle d_{i}\rangle_{\alpha}\in(0,1]$ we have
\[ \lim_{n\to\infty} \frac{1}{n}(d_{1}+d_{2}\dots +d_{n})=\frac{\alpha}{\alpha-1}\]
and
\[ \lim_{n\to\infty} \sqrt[n]{d_{1}d_{2}\dots d_{n}}=\prod_{d=1}^\infty \sqrt[\alpha^d]{d^{\alpha-1}}.\]
\end{theorem}
\proof
Let 
\[ f(x)=\sum_{d=1}^{\infty}d\chi_{\left(\frac{1}{\alpha^{d}},\frac{1}{\alpha^{d-1}}\right]}(x).\] 
By lemma 3.1 we have $f(T^{i-1}(x))=d_{i}$ for all $i\in\mathbb{N}$, where $x=\langle d_{i}\rangle_{\alpha}$.
Applying Birkhoff's ergodic theorem we obtain
\[ \lim_{n\to\infty}  \frac{1}{n}\sum_{i=1}^{n}d_{i}=\lim_{n\to\infty}\frac{1}{n}\sum_{i=1}^{n}f(T^{i-1}(x))=\int_{0}^{1} f(x)dx \]
\[ = \sum_{d=1}^{\infty}d\frac{\alpha-1}{\alpha^{d}}=\frac{\alpha}{\alpha-1}\]
for almost all $x=\langle d_{i}\rangle_{\alpha}\in(0,1]$. Now let 
\[ f(x)=\sum_{d=1}^{\infty}\log(d)\chi_{\left(\frac{1}{\alpha^{d}},\frac{1}{\alpha^{d-1}}\right]}(x).\] 
By lemma 3.1 we here have $f(T^{i-1}(x))=\log(d_{i})$ for all $i\in\mathbb{N}$. Applying Birkhoff's ergodic theorem again we obtain
\[ \lim_{n\to\infty}  \frac{1}{n}\sum_{i=1}^{n}\log(d_{i})=\lim_{n\to\infty}\frac{1}{n}\sum_{i=1}^{n}f(T^{i-1}(x))=\int_{0}^{1} f(x)dx \]
\[ = \sum_{d=1}^{\infty}\log(d)\frac{\alpha-1}{\alpha^{d}}\]
for almost all $x=\langle d_{i}\rangle_{\alpha}\in(0,1]$. Applying the exponential on both sides gives the result on the geometric mean.\qed~\\~\\
Consider the function
\[ G(x)=\prod_{d=1}^\infty \sqrt[\alpha^x]{d^{x-1}}.\]
$G(2)= 1.6616879496\dots$ is known as  Somos' constant, see page 446 of \cite{[FI]}. In definition 2 of \cite{[SH]} the numbers $\sqrt[x-1]{G(x)}$ are called generalized Somos' constants. We would prefer to call $G:(0,\infty)\to \mathbb{R}$ the Somos' function. We display an approximation of this function in figure 1.  We like to remark that from theorem 8 in \cite{[SH]} we get a relationship between $G$ and the generalized-Euler-constant function $\gamma$:
\[G(x)=\frac{x}{x-1}e^{-\gamma(1/x)/x},\]
where 
\[ \gamma(x)=\sum_{i=1}^{\infty}x^{i-1}\left(\frac{1}{i}-\log(\frac{i+1}{i})\right).\]      
The proof of the formula is just a straightforward calculation. 
\begin{figure}
\vspace{0pt}\hspace{0pt}\scalebox{0.5}{\includegraphics{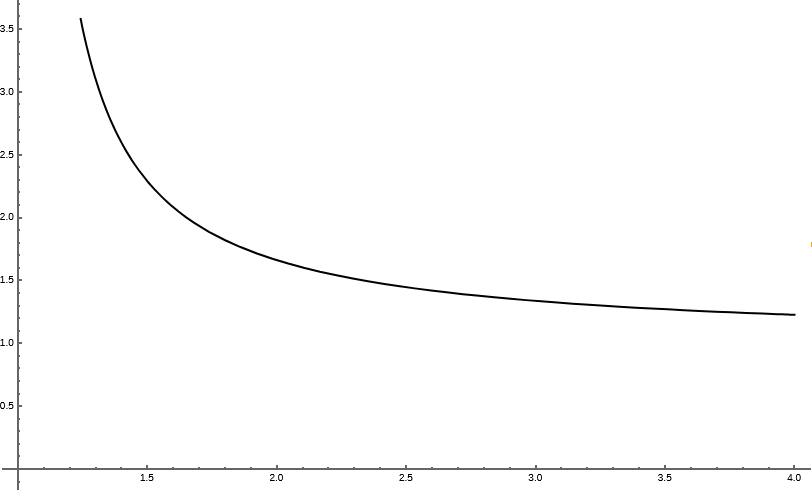}}
\caption{The Somos' function $G$}
\end{figure}
\section{Dimension theoretical aspects}
For a proper subset $D\subset\mathbb{N}$ and $\alpha>1$ let us consider the set of numbers which have only digits in $D$ in their $\alpha$-expansion,
\[ I_{\alpha}(D)=\{\langle d_{i}\rangle_{\alpha}|d_{i}\in D\}.\] 
This is a Cantor set hence its Hausdorff dimension $\dim_{H}I_{\alpha}(D)$ is of interest. The reader who is not familiar with basic dimension theory should consider our essential \cite{[NE2]} or \cite{[FA]}. Applying well known results in dimension theory we obtain:
\begin{theorem}
Let $D\subset\mathbb{N}$, $\alpha>1$ and $h\ge 0$ be the unique solution of
\[ \sum_{i\in D}(\frac{\alpha-1}{\alpha^{i}})^{h}=1.\]
Than $\dim_{H}I_{\alpha}(D)=h$.
\end{theorem}
\proof Let $\tilde I_{\alpha}(D)=I_{\alpha}(D)\cup\{0,0\}$, adding a point does not change Hausdorff dimension. $\tilde I_{\alpha}(D)$ is the attractor of the iterated function system $([0,1],\{T_{i}|i\in D)$, since
\[ \tilde I_{\alpha}(D)=\bigcup_{i\in D}  T_{i}(\tilde I_{\alpha}(D)).\]
This iterated function system fulfills the open set condition since
\[ T_{i}((0,1))\cap T_{j}((0,1))=\emptyset \]
if $i\not=j$. The dimension formula in our result is the classical Moran formula. If $A$ is finite the result directly follows from the work of Moran \cite{[MO]}. If $A$ is infinite it follows from theory of infinite iterated function systems see theorem 3.11 of \cite{[HF]} or \cite{[MU]}.\qed~\\~\\
For $D_{n}=\{1,\dots,n\}$ we obtain $\dim_{H}I_{\alpha}(D_{n})=h(n)$ where $h(n)>0$ is given by
\[ (1-\alpha^{-h(n)n})\frac{(\alpha-1)^{h(n)}}{\alpha^{h(n)}-1}=1.\]
Note that $h(n)\to\infty$ for $n\to\infty$ hence
\[ \dim_{H}\{\langle d_{i}\rangle_{\alpha}|(d_{i})\mbox{ is bounded}\}=\dim_{H}\bigcup_{n=1}^{\infty}I_{\alpha}(D_{n})=1,\]
although the Lebesgue measure of this set is zero.~\\~\\ 
Next we study the set of real numbers with prescribed frequencies of digits in their $\alpha$-expansion. We need some notations. Let $\mathfrak{p}=(p_{1},\dots,p_{n})\in[0,1)^n$ be a probability distribution on $\{1,\dots,n\}$ that means
\[ \sum_{i=1}^{n}p_{i}=1.\]
The expected value of $\mathfrak{p}$ is
\[ E(\mathfrak{p})=\sum_{i=1}^{n}ip_{i}\]
and the entropy is given by
\[ H(\mathfrak{p})=-\sum_{i=1}^{n}p_{i}\log(p_{i}).\]
We are interested in the set of real numbers with frequency $p_{d}$ of the digit $d$ in their $\alpha$-expansion;
\[ \mathfrak{F}_{\alpha}(\mathfrak{p})=\{x\in(0,1]|\mathfrak{f}_{d,\alpha}(x)=p_{d},~d=1,\dots, n \}.\]   
These are fractal sets and we will proof the following dimension formula:
\begin{theorem}
For all $\alpha>1$ and all probability distribution $\mathfrak{p}\in[0,1)^n$ on $\{1,\dots,n\}$ we have
\[ \dim_{H}\mathfrak{F}_{\alpha}(\mathfrak{p})=\frac{ H(\mathfrak{p})}{\log(\alpha)E(\mathfrak{p})-\log(\alpha-1)}.\]
\end{theorem}
\proof Let $P$ be the probability measure on $\mathbb{N}^\mathbb{N}$ which is the product of the probability distribution $\mathfrak{p}=(p_{1},\dots,p_{n},0,0,\dots)$ on $\mathbb{N}$.
We project $P$ to a probability measure $\mu$ on $(0,1]$ using $\pi$,
\[ \mu=\pi(P)=P\circ \pi^{-1}.\]
Be the law of large numbers for $P$-almost all sequences $(d_{i})$ the frequency of a digit $d$ is given by $p_{d}$, hence
$\mu(\mathfrak{F}_{\alpha}(\mathfrak{p}))=1$.~\\
The intervals
\[ I_{d_{1},\dots,d_{k}}=T_{d_{k}}\circ \dots\circ T_{d_{1}}((0,1])\]
have length and measure given by
\[ |I_{d_{1},\dots,d_{k}}|=\frac{(\alpha-1)^{k}}{ \alpha^{d_{1}+\dots+d_{k}}},\quad \mu(I_{d_{k},\dots,d_{1}})=\prod_{j=1}^{k}p_{d_{j}},\]
compare with section 2.1. For $x=\langle d_{i}\rangle_{\alpha}\in(0,1]$ let $I_{d_{1},\dots,d_{k}}(x)$ the interval containing $x$ and
\[ \sharp_{d}(x|k)=\sharp\{j|d_{j}=d,\quad j=1,\dots, k\}.\] 
For the length of the intervals we obtain 
\[ \log(|I_{d_{1}\dots d_{k}}(x)|)=\log(\frac{(\alpha-1)^{k}}{ \alpha^{d_{1}+\dots+d_{k}}})=k\log(\alpha-1)-\log(\alpha)\sum_{d=1}^{n}\sharp_{d}(x|k)d\]
and for the measure of the intervals we get
\[ \log(\mu(I_{d_{1}\dots d_{k}})(x))=\log(\prod_{j=1}^{k}p_{d_{j}})=\sum_{d=1}^{n} \sharp_{d}(x|k)\log p_{d}.\]
If $x\in\mathfrak{F}_{\alpha}(\mathfrak{p})$, we have
\[\lim_{k\to\infty} \frac{\sharp_{d}(x|k)}{k}=p_{d}\]
for all $d\in\mathbb{N}$. Hence
\[ \lim_{k\longmapsto \infty} \frac{1}{k}\log \frac{\mu(I_{n_{1}\dots n_{k}}(x))}{|I_{n_{1}\dots n_{k}}(x)|^s}\]
\[=-H(\mathfrak{p})+s(\log(\alpha)E(\mathfrak{p})-\log(\alpha-1)).\]
This implies 
\[ \lim_{k\longmapsto \infty}  \frac{\mu(I_{n_{1}\dots n_{k}}(x))}{|I_{n_{1}\dots n_{k}}(x)|^s}=\{\begin{array}{cc} 0& s<h\\ \infty & s>h \end{array},\]
where
\[ h=\frac{ H(\mathfrak{p})}{\log(\alpha)E(\mathfrak{p})-\log(\alpha-1)}.\]
By the local mass distribution principle, see proposition 4.9 of \cite{[FA]}, we have $\mathfrak{H}^{s}(\mathfrak{F}_{\alpha}(\mathfrak{p}))=\infty$ for $s<h$ and $\mathfrak{H}^{s}(\mathfrak{F}_{\alpha}(\mathfrak{p}))=0$ for $s>h$, where $\mathfrak{H}^{s}$ is the $s$-dimensional Hausdorff measure. This implies
$\dim_{H}(\mathfrak{F}_{\alpha}(\mathfrak{p}))=h$.
  \qed
\section{The base-change-transformation} 
For fixed real numbers $\alpha,\beta>1$ we define a base-change-transformation by $f:[0,1]\to[0,1]$ by $f(x)=f(\langle d_{i}\rangle_{\alpha})=\langle d_{i}\rangle_{\beta}$, that is
\[ f(x)=f(\sum_{i=1}^{\infty}(\alpha-1)^{i-1}\alpha^{-(d_{1}+\dots+d_{i})})=\sum_{i=1}^{\infty}(\beta-1)^{i-1}\beta^{-(d_{1}+\dots+d_{i})}\]
for $x\not =0$ and $f(0)=0$. Obviously we have the functional equation
\[ f(x/\alpha)=f(x)/\beta \]
for all $x\in[0,1]$. If $\alpha=\beta$ the map $f$ is the identity on $[0,1]$. In general we have: 
\begin{theorem}
$f:[0,1]\to[0,1]$ is strictly increasing and continuous.
\end{theorem}
\proof Let us introduce some notations. For $(d_{i})\in\mathbb{N}^{\mathbb{N}}$ we set
\[\Xi(d_{i})=\min\{k|d_{i}=1 \forall i>k\}\]
if the minimum exists and $\Xi(d_{i})=\infty$ if this is not the case. For two different sequences  $(d_{i}),(g_{i})\in\mathbb{N}^{\mathbb{N}}$ we set
\[ |(d_{i})\wedge (g_{i})|=\min\{i|d_{i}\not= g_{i}\}.\]   
Let $x,y\in (0,1]$ with $x=\langle d_{i}\rangle_{\alpha}$ and $y=\langle g_{i}\rangle_{\alpha}$. If $x<y$, we have \
$d_{|(d_{i})\wedge (g_{i})|}>g_{|(d_{i})\wedge (g_{i})|}$. Hence $f(x)=\langle d_{i}\rangle_{\beta}<\langle g_{i}\rangle_{\beta}=f(y)$. If $x>0$ obviously $f(x)>f(0)=0$. Hence $f$ is strictly increasing. \\
Let $(x_{n})=\langle d_{i}(n)\rangle_{\alpha} \in (0,1]^{\mathbb{N}}$ be a sequence. If $\lim_{n\to\infty} x_{n}=0$, we have $\lim_{n\to\infty} d_{i}(x_{n})=\infty$ for all $i\in\mathbb{N}$. Hence 
\[ \lim_{n\to\infty}f(x_{n})=\lim_{n\to\infty}\langle d_{i}(n)\rangle_{\beta}=0=f(0)\]
and $f$ is continuous in $x=0$. \\
Now let $x=\langle d_{i}\rangle_{\alpha}\in(0,1]$ and $\lim_{n\to\infty}x_{n}=x$. If  $\Xi(d_{i})=\infty$, 
we have $\lim_{n\to\infty}|(d_{i}(n))\wedge (d_{i})|=\infty$ hence
\[ \lim_{n\to\infty}f(x_{n})=\lim_{n\to\infty}\langle d_{i}(n)\rangle_{\beta}=\langle d_{i}\rangle_{\beta}=f(x).\]
Now let $\Xi(d_{i})=k$ and let $(n_{j})$ be a sequence of natural numbers with $\lim_{k\to\infty}n_{j}=0$. There are to possibilities.
Either $\lim_{j\to\infty}|(d_{i}(n_{j}))\wedge (d_{i})|=\infty$ or $d_{k}(n_{j})=d_{k}-1$ for $j$ sufficient large and $\lim_{j\to\infty}d_{i}(n_{j})=\infty$ for all $i>k$. In both cases
\[ \lim_{j\to\infty}f(x_{n_{j}})=\lim_{j\to\infty}\langle d_{i}(n_{j})\rangle_{\beta}=\langle d_{i}\rangle_{\beta}=f(x).\]
Hence $f$ is continuous for all $x\in(0,1]$.
\qed~\\\\
Now we have a look on the structure of the graph of $f$, \[ F=\{(x,f(x))|x\in[0,1]\}\subseteq [0,1]\times [0,1],\]
\[ =\{0,0\}\cup\{\sum_{i=1}^{\infty}(\alpha-1)^{i-1}\alpha^{-(d_{1}+\dots+d_{i})},\sum_{i=1}^{\infty}(\beta-1)^{-1}\beta^{-(d_{1}+\dots+d_{i})})|(d_{i})\in\mathbb{N}^{\mathbb{N}}\}.\] We display an approximation of $F$ for $(\alpha,\beta)=(3,2),(3,10)$ in figure 2.
\begin{figure}
\vspace{0pt}\hspace{0pt}\scalebox{0.35}{\includegraphics{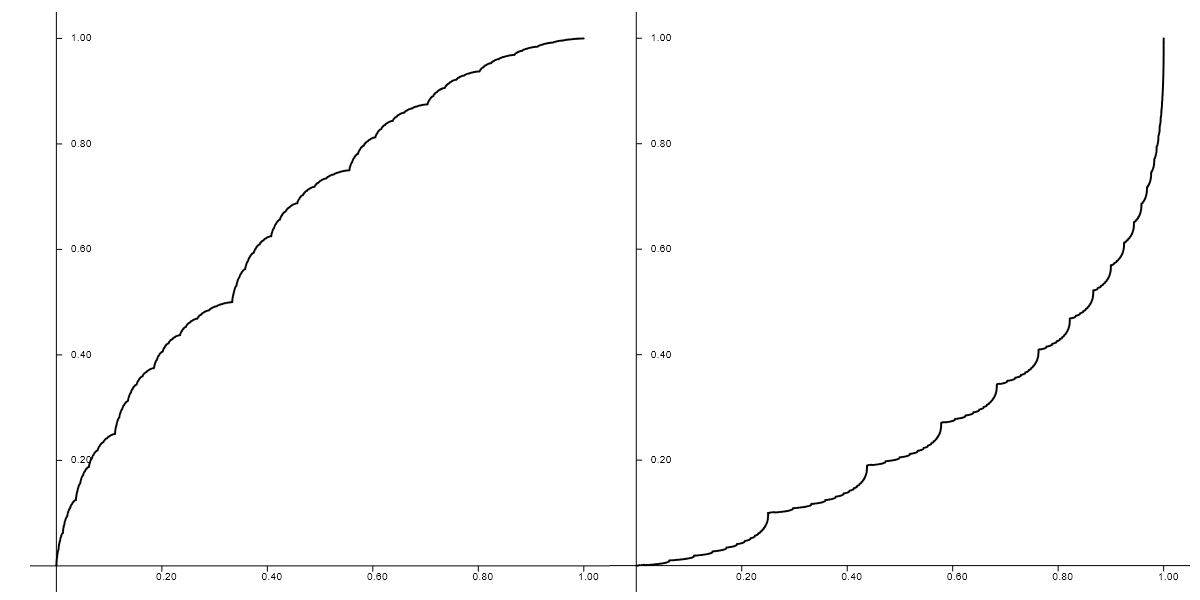}}
\caption{An approximation of the base-change-transformation for $(\alpha,\beta)=(3,2),(3,10)$. }
\end{figure}
In addition to theorem 5.1 we obtain:
\begin{theorem}
$F$ is an infinite self-affine curve with Hausdorff dimension one.
\end{theorem}
\proof
Consider affine contractions $G_{i}:[0,1]\times [0,1]\to[0,1]\times [0,1]$, given by 
\[ G_{i}(x,y)=(\frac{\alpha-1}{\alpha^{i}}x+\frac{1}{\alpha^{i}}, \frac{\beta-1}{\beta^{i}}y+\frac{1}{\beta^{i}})\]
for $i\in\mathbb{N}$. We have
\[ \bigcup_{i=1}^{\infty}G_{i}(F)=F,\]
which means that $F$ is infinite self-affine. We have $\dim_{H}F\ge 1$ since the projection of $F$ to the first coordinate axis is the whole interval and 
projections do not increase Hausdorff dimension. Now consider the finite iterated function system, given by $([0,1]^2,\{G_{1},G_{\infty}\})$, where $G_{\infty}$ is given by 
\[  G_{\infty}(x,y)=(\frac{1}{\alpha}x,\frac{1}{\beta}y).\]
Since $G_{i}([0,1]^2)\subseteq G_{\infty}([0,1]^2)$ for all $i\ge 2$ the attractor $\Lambda$ of this finit iterated function system contains $F$. But it follows from the theory of finit self-affine sets, that $\dim_{H}\Lambda=1$, see proposition 10.2.6 in \cite{[BSS]}. Hence $\dim_{H}F=1$.\qed~~\\\\
At the end we like to include a result on the upper derivative $\overline{D}_{x}f$ resp. lower derivative $\underline{D}_{x}f$ of the base-change-transformation $f$.
\begin{theorem}
For $\alpha,\beta>1$ with $\alpha\not=\beta$ let $f$ be the corresponding base-change-transformation. 
For almost all $x\in(0,1]$ we have $\overline{D}_{x}f=\infty$ and $\underline{D}_{x}f=0$.
\end{theorem} 
\proof 
Consider $f(x)=\langle d_{i}\rangle_{\beta}\in(0,1]$ with
\[ \lim_{n\to\infty} \frac{1}{n}(d_{1}+d_{2}\dots +d_{n})=\frac{\beta}{\beta-1}.\]
From theorem 3.2. we know that the set of $x=\langle d_{i}\rangle_{\alpha}\in(0,1]$ with these property has full measure. Define
\[ x_{n}=T_{d_{n}}\circ \dots\circ T_{d_{1}}(0)\mbox{ and }y_{n}=T_{d_{n}}\circ \dots\circ T_{d_{1}}(x)(1)\]
We obviously have $x_{n}<x<y_{n}$ and $\lim_{n\to\infty}x_{n}=x$ as well as $\lim_{n\to\infty}y_{n}=x$. Furthermore we obtain
\[ \lim_{n\to\infty}\frac{f(y_{n})-f(x_{n})}{y_{n}-x_{n}}=\lim_{n\to\infty}\frac{(\beta-1)^{n}\beta^{-d_{1}-\dots-d_{n}}}{(\alpha-1)^{n}\alpha^{-d_{1}-\dots-d_{n}}}
=\lim_{n\to\infty}\left(\frac{(\beta-1)\alpha^{(d_{1}+\dots+d_{n})/n}}{(\alpha-1)\beta^{(d_{1}+\dots+d_{n})/n}}\right)^{n}=\infty.\]
since
\[\frac{(\beta-1)\alpha^{\beta/(\beta-1)}}{(\alpha-1)\beta^{\beta/(\beta-1)}}>1\] 
for $\alpha,\beta>1$ with $\alpha\not=\beta$. Hence $\overline{D}_{x}f=\infty$ for almost all $x\in(0,1]$. Furthermore this implies $\underline{D}_{x}f^{-1}=0$ for almost all $x\in(0,1]$, where $f^{-1}$ is the base-change-transformation from $\beta$ to $\alpha$. Interchanging $\alpha$ and $\beta$ the stated result follows for $(\alpha,\beta)$ and $(\beta,\alpha)$ since the intersection of sets with full measure has full measure.

\qed


\begin{thebibliography}
\small
\bibitem{[DY]} D. Airey and B. Mance, 
The Hausdorff dimension of sets of numbers defined by their Q-Cantor series expansions,
J. Fractal Geom. 3, No. 2, 163-186, 2016.
\bibitem{[BA]} S. Baker, The growth rate and dimension theory of beta-expansions, Fundam. Math. 219, No. 3, 271-285, 2012.
\bibitem{[BSS]} B. Barany, K. Simon and B. Solomyak,
Self-similar and self-affine sets and measures, Mathematical Surveys and Monographs 276. Providence, 2023. 
\bibitem{[BAR]} L. Barreira and F. Iommi, 
Frequency of digits in the Lüroth expansion, J. Number Theory 129, No. 6, 1479-1490, 2009. 
\bibitem{[BE]} A.S. Besicovitch,  On the sum of digits of real numbers represented in the dyadic system, Math. Ann. 110, 321-330, 1934.
\bibitem{[BO]} E. Borel, Les probabilités denombrables et leurs applications arithmetiques, Rend. Circ. Mat. Palermo 27, 247–271, 1909.
\bibitem{[EIN]} M. Einsiedler and T. Ward,  Ergodic theory. With a view towards number theory, Graduate Texts in Mathematics 259. Springer London, 2011. 
\bibitem{[EG]} H.G. Eggleston, {\it The fractional dimension of a set defined by decimal properties},
Quart. J. Math. Oxford Ser. 20, 6-31, 1949.
\bibitem{[ER]} P. Erdos and A. Rényi, 
Some further statistical properties of the digits in Cantor’s series, Acta Math. Acad. Sci. Hung. 10, 21-29, 1959. 
\bibitem{[GA]} Galambos J, Representations of Real Numbers by Infinite Series, Lecture Notes in Math. 502, Springer, 1976.
\bibitem{[JA]} V. Jarnik,  Zur metrischen Theorie der diophantischen Approximationen, Prace Mat.-. Fix. 36, 91-106, 1928. 
\bibitem{[FANG]} L. Fang and M. Wu, Hausdorff dimension of certain sets arising in Engel expansions, Nonlinearity 31, 2105–2125, 2018.
\bibitem{[FA]} K. Falconer, Fractal Geometry - Mathematical Foundations and Applications, Wiley, New York, 1990.
\bibitem{[HF]} H. Fernau, Infinite iterated function systems, Math. Nach., vol. 170, issue 1, 79-91, 1994.
\bibitem{[FI]} S. R. Finch, Mathematical Constants, Cambridge University Press, Cambridge, 2003.
\bibitem{[KIF]} Y. Kifer Fractal dimensions and random transformations, Trans. Amer. Math. Soc. 348, no. 5, 2003-2038, 1996.
\bibitem{[KI]} A. Khinchin, Metrische Kettenbruchprobleme, Compositio Mathematica 1 (1935) 361-382.
\bibitem{[KO]} V. Komornik, J. Lu and Y. Zou, Expansions in multiple bases over general alphabets, Acta Math. Hung. 166, No. 2, 481-506, 2022. 
\bibitem{[HA]} F. Hausdorff, Dimension und äußeres Maß, Mathematische Annalen. 79 (1–2), 157–179, 1919. 
\bibitem{[JAG]} H. Jager H. and C. de Vroedt, L\"uroth series and their ergodic properties, Indag. Math. 31, 31–42, 1969.
\bibitem{[LY]} A. Lasota and J.A.Yorke,  On the existence of invariant measures for piecewise monotonic transformations,
 Transactions of the American Mathematical Society 186, 481-488, 1973.
\bibitem{[LI]} Y. Li, Expansions in multiple bases, Acta. Math. Hungar., 163, 576-600, 2021.
\bibitem{[MO]} P. Moran, Additive functions of intervals and Hausdorff measure, Proc. Camb. Phil. Soc. 42, 15-23, 1946.
\bibitem{[MU]} R.D. Mauldin and M. Urbanski, Dimensions and measures in infinite iterated function systems, Proc. London Math. Soc., 73, 105-154, 1996.
\bibitem{[NE3]} J. Neunh\"auserer, Non-uniform expansions of real numbers, Mediterr. J. Math., 18, 8pp., 2021.
\bibitem{[NE4]} J. Neunh\"auserer, Expansions of real numbers with respect to two integer bases, arXiv:2002.10824, 2020; to appear Rocky Mountain Journal of Mathematics, 2024. 
\bibitem{[NE1]} J. Neunh\"auserer, Introduction to ergodic theory, Essentials, Springer Spektrum, Wiesbaden, 2020.
\bibitem{[NE2]} J. Neunh\"auserer, Dimensional theory, Essentials, Springer Spektrum, Wiesbaden, 2021.
\bibitem{[RE]}  A. Renyi, Representations for real numbers and their ergodic properties. Acta Math.
Acad. Sci. Hungar 8, 477–493, 1957. 
\bibitem{[SH]} J. Sondow and P. Hadjicostas, The generalized Euler-constant function $\gamma(z)$ and a generalization of Somos's quadratic recurrence
constant, J. Math. Anal. Appl. 332, 292-314, 2007.
\bibitem{[WU]} J. Wu. A problem of Galambos on Engel expansions, Acta Arith. 92, No. 4, 383-386, 2000. 
\end{thebibliography}
\end{document}